\documentclass
[12pt,a4paper]{article}
\usepackage[cp1251]{inputenc}
\usepackage{amssymb, amsfonts, amsmath, latexsym}
\usepackage[english]{babel}

\begin{document}

\normalsize\centerline{\bf  I. V.  Protasov }\vspace{6 mm}

\Large
\centerline{\bf  A note on bornologies }\vspace{6 mm}
\normalsize

\centerline{\it To Taras Banakh on  50th birthday}
\vspace{4 mm}

{\bf Abstract.} A bornology on a set $X$ is a family   $\mathcal{B}$ of subsets of $X$  closed under taking subsets, finite unions and such that $\cup \mathcal{B}=X$.

We prove that, for a bornology $\mathcal{B}$ on $X$, the following statements are equivalent:

(1)	there exists a vector topology  $\tau$ on  the vector space  $\mathbb{V} (X) $  over $\mathbb{R}$
 such that $\mathcal{B}$ is the family of all subsets of $X$ bounded in $\tau$;

(2)	there exists a uniformity $\mathcal{U}$   on $X$  such  that $\mathcal{B}$ is the family of all subsets of $X$ totally bounded in $\mathcal{U}$;

(3)	for every  $Y \subseteq X$,  $Y \notin \mathcal{B}$, there exists a metric $d$ on $X$  such that $\mathcal{B}\subseteq \mathcal{B} _{d}$,  $Y\notin  \mathcal{B} _{d}$, where $\mathcal{B} _{d}$  is the family of all closed discrete subsets of $(X, d)$;

(4)	for every     $Y \subseteq X$,  $Y \notin \mathcal{B}$, there exists
 $Z\subseteq  Y$  such that
  $Z^{\prime} \notin \mathcal{B}$  for each infinite subset $Z^{\prime}$  of $Z$.

A bornology $\mathcal{B}$  satisfying $(4)$  is called antitall.  We give topological and functional characterizations of antitall bornologies.
\vspace{6 mm}

UDC 519.51

\vspace{3 mm}


\vspace{6 mm}

A family $\mathcal{I}$ of subsets of a set $X$  is called an {\it ideal} (in the Boolean algebra $\mathcal{P}_{X}$ of all subsets of $X$)  if $\mathcal{I}$  is closed under formations of finite unions and subsets.
If $\cup \mathcal{I}=X$  then $\mathcal{I}$ is called a {\it bornology}, so a bornology is an ideal containing the ideal $\mathcal{F}_{X}$ of all finite subsets of  $X$.

For an ideal  $\mathcal{I}$, a family $\mathcal{F}\subseteq \mathcal{I}$ is called a {\it base} of $\mathcal{I}$  if, for any $A\in \mathcal{I}$, there exists $B\in \mathcal{F}$  such that $A\subseteq B$.

An ideal $\mathcal{I}$ on $X$ is called {\it tall} if any infinite subset $Y$  of $X$  contains an infinite subset $Z$  such that $Z\in \mathcal{I}$.

We say that an ideal $\mathcal{I}$ on $X$  is {\it antitall} if, for any $Y\subseteq  X$,  $Y\notin \mathcal{I}$ there exists $Z\subseteq Y$  such that $Z^{\prime}\notin \mathcal{I}$  for each infinite subset $Z^{\prime}$ of $Z$.
Clearly, every bornology with a countable base is antitall.
In particular, a bornology of all bounded subsets of a metric space is antitall.
On the  other hand, for every bornology $\mathcal{B}$
 with countable base,  there is a metric $d$ on $X$ such that $\mathcal{B}$ is the bornology of bounded subsets of
 $(X, d)$.

Every   bornology is the meet of some tall and antitall bornologies,  see Proposition 1.

Given a  bornology $\mathcal{B}$  on $X$ and
a set $\mathcal{S}$    of bornologies on $X$, we say that  $\mathcal{B}$  is {\it approximated} by $\mathcal{S}$ if,
for every $Y\subseteq  X$,  $Y\notin  \mathcal{B}$,  there exists  $\mathcal{B}^{\prime}\in  \mathcal{S}$ such that
$\mathcal{B}\subseteq  \mathcal{B}^{\prime}$  and $Y\notin \mathcal{B}^{\prime}$.
If $\mathcal{B}$  is approximated  by $\mathcal{S}$ and  $\cap \mathfrak{F }\in  \mathcal{S}$ for any $\mathfrak{F}\subseteq  \mathcal{S}$  then $\mathcal{B}\in \mathcal{S}$.

 We use the standard topological terminology, see [1].
For the history of bornology, see [2].

\section{Vector  bornologies}

For a set $X$,  $\mathbb{V(}X) $  denotes the vector space over  $\mathbb{R}$  with the basis $X$.
Under a {\it vector topology}  on $\mathbb{V}(X) $,  we mean a topology $\tau$ such that
$(\mathbb{V}(X), \tau )$
 is a topological vector space.

Let $\tau$  be a vector topology  on $X$,  $U$ be a neighbourhood  of $0$  in  $\tau$.
A subset $Y\subseteq  X$  is called  $U$-{\it bounded}  if there exists $n\in \mathbb{N}$ such that $Y\subseteq n U$.
If $Y$ is  $U$-bounded for every neighbourhood $U$  of $0$  in $\tau$  then $Y$  is called
$\tau$-bounded.
We denote by $\mathcal{B}_{\tau}$  the family of all $\tau$-bounded  subsets of  $X$  and observe that
$\mathcal{B}_{\tau}$
 is a   bornology.
\vspace{3 mm}

{\bf Theorem 1. } {\it For a bornology $\mathcal{B}$  on a set $X$,  the following statements are equivalent:
\vspace{2 mm}

$(i)$  $\mathcal{B}= \mathcal{B}_{\mu}$   for some vector topology $\mu$  on $\mathbb{V}(X)$;\vspace{2 mm}

$(ii)$  $\mathcal{B}$  is approximated by the set $\{\mathcal{B}_{\tau}: \tau$  is a vector topology on
$\mathbb{V}(X)$;\vspace{2 mm}

$(iii)$  $\mathcal{B}$ is antitall.
\vspace{5 mm}

Proof.}
The implication $(i)\Longrightarrow  (ii)$ is evident.
To see $(ii)\Longrightarrow  (i)$, we denote by $\mu$ the strongest vector topology
 such that each subset $Y\in \mathcal{B}$  is $\mu$-bounded.
\vspace{5 mm}

$(ii)\Longrightarrow  (iii)$ Let $Y\subseteq X$ and $Y\notin \mathcal{B}$.  We choose a vector  topology $\tau$
such that $\mathcal{B}\subseteq  \mathcal{B}_{\tau}$, $Y\notin\mathcal{B}_{\tau}$.
Then there exists a neighbourhood $U$  of  $0$ in $\tau$ such that
 $Y \setminus n \ U \neq\emptyset $  for
  each $n\in \mathbb{N}$.
We take a subset
$Z= \{z_{n} : n\in \mathbb{N}\}$ of $Y$  such that $z_{n}\notin n \  U$.
Then $Z^{\prime}\notin \mathcal{B}_{\tau}$  for each infinite subset  $Z^{\prime}$ of $Z$.
Since
$\mathcal{B}\subseteq\mathcal{B}_{\tau}$, we conclude that $\mathcal{B}$ is antitall.

\vspace{5 mm}

$(iii)\Longrightarrow  (ii)$.
We write $\mathbb{V}(X)$  as  $\oplus \mathbb{R} _{x}, x\in X$, $\mathbb{R} _{x}$ is a copy of $\mathbb{R}$,
  so each vector $v\in \mathbb{V}(X)$ is of the form
  $(\lambda_{x})_{x\in X}$, $\lambda_{x}\in \mathbb{R}_{x}$
    and $\lambda_{x}=0$ for all but finitely many $x\in X$.

Let $Y \subseteq X$  and $Y\notin \mathcal{B}$. We take a countable subset $Z$ of $Y$  such that $Z^{\prime}\notin \mathcal{B}$  for each infinite subset $Z^{\prime}$  of $Z$.
We denote by  $\mathcal{P}$  the partition of $X$ into subsets
$X\setminus Z$,  $\{z\}$, $z\in Z$,  and by
$\Lambda_{\mathcal{P}}$ the family of all functions
 $\lambda: \mathcal{P}\longrightarrow\{\frac{1}{n}: n\in\mathbb{N}\}$.
For each $\lambda\in\Lambda_{\mathcal{P}}$, we put
$$U(\lambda)=\{v\in \mathbb{V}(X): v=(\lambda_{x})_{x\in X},   \   \   |\lambda_{x}|< \lambda(P)  \  for  \  each \  x\in P, \  P\in\mathcal{P}\}, $$
 and take
 $\{U(\lambda): \lambda\in\Lambda_{\mathcal{P}}\}$
  as a base at $0$  for some (uniquely determined) vector topology $\tau$  on  $\mathbb{V}(X) $.
By the construction, a subset $A$ of $X$  is  $\tau$-bounded if and only if  $A=\cup \mathfrak{F}$  for some finite subset  $\mathfrak{F}$ of $\mathcal{P}$.
It follows that $\mathcal{B}\subseteq\mathcal{B}_{\tau}$ and  $Z\notin  \mathcal{B}_{\tau}$ so $Y\notin\mathcal{B}_{\tau}$.   $ \ \  \  \Box$

\section{Totally bounded  bornologies}

For a uniformity  $\mathcal{U}$  on a set  $X$,  $x\in X$  and an entourage $\varepsilon\in\mathcal{U}$,
 the set $B(x,\varepsilon)=\{y\in X: (x,y)\in\varepsilon\}$  is called an
 $\varepsilon$-{\it ball} centered at $x$.
A subset $Y\subseteq X$  is called {\it totally bounded}  in the uniform space $(X, \mathcal{U}) $  if, for each
 $\varepsilon\in \mathcal{U}$,  $Y$  can be covered by finite number  of
  $\varepsilon$-balls.
If $Y$ is not totally bounded then there exists an  $\varepsilon$-discrete subset
$\{y_{n}: n\in \omega\}$ of $X$,  i.e.  $B(y_{n}, \varepsilon)\cap  B(y_{m}, \varepsilon)= \emptyset$
           for all distinct $n, m\in \omega$.

Given a uniformity $\mathcal{U}$ on $X$,  we denote by $\mathcal{B}_{\mathcal{U}}$ the
 bornology of all totally  bounded subsets of $(X, \mathcal{U})$.
\vspace{4 mm}

{\bf Theorem 2.} {\it For a bornology $\mathcal{U}$  on a set $X$, the following  statements are equivalent:
\vspace{2 mm}

$(i)$  $\mathcal{B}= \mathcal{B}_{\mathcal{M}}$  for some uniformity $\mathcal{M}$ on $X$;\vspace{2 mm}

$(ii)$  $\mathcal{B}$ is approximated by the set $\{ \mathcal{B}_{\mathcal{U}}: \mathcal{U} $  is a uniformity on $X$\};\vspace{2 mm}

$(iii)$  $\mathcal{B}$ is antitall.

\vspace{5 mm}
Proof.}
The implication $(i)\Longrightarrow (ii)$  is evident.
To see $(ii)\Longrightarrow (i)$,  we take the strongest uniformity $\mathcal{M}$ such that each
subset $Y\in \mathcal{B}$   is totally bounded in $(X, \mathcal{M})$.
\vspace{3 mm}

$(ii)\Longrightarrow (iii)$ .
Let $Y\subseteq X$  and $Y\notin \mathcal{B}$.
We choose a uniformity $\mathcal{U}$ on $X$  such that $\mathcal{B} \subseteq \mathcal{B}_{\mathcal{U}}$  and $Y\notin \mathcal{B}_{\mathcal{U}}$.
Then there exists  $\varepsilon\in \mathcal{U}$ such that some countable subset $Z$ of $Y$ is $\varepsilon$-discrete.  Clearly, $Z^{\prime}\notin \mathcal{B}_{\mathcal{U}}$ for each infinite subset $Z^{\prime}$  of $Z$. Since $\mathcal{B}\subseteq \mathcal{B}_{\mathcal{U}}$, we conclude that $\mathcal{B}$ is antitall.
\vspace{3 mm}

$(iii)\Longrightarrow (i)$.
Let $Y\subseteq X$  and $Y\notin \mathcal{B}$.
We take a countable subset $Z$ of $Y$  such that $Z^{\prime}\notin  \mathcal{B}$  for each countable subset  $Z^{\prime}$ of $Z$.
Then we define a uniformity $\mathcal{U}$  on $X$  such that
 $X\setminus  Z$  is totally bounded in $\mathcal{U}$  and $Z$  is $\varepsilon$-discrete for some $\varepsilon\in \mathcal{U}$.
Clearly, $\mathcal{B}\subseteq \mathcal{B}_{\mathcal{U}}$  and $Y\notin \mathcal{B}_{\mathcal{U}}$. $ \ \ \ \  \Box$

\section{Closed discrete   bornologies}

For a metric space $(X, d)$, we denote by $\mathcal{B}_{d}$ the bornology of all closed discrete  subsets of $X$.
Clearly, $Y\notin \mathcal{B}_{d}$ if and only if there exists an injective sequence in $Y$  converging in $X$.
\vspace{3 mm}

{\bf Theorem 3. } {\it For every bornology $\mathcal{B}$  on a set $X$,  the following statements are equivalent:
\vspace{2 mm}

$(i)$  $\mathcal{B}$  is approximated by the set $\{\mathcal{B}_{d}: d $  is a metric on $X\}$;
\vspace{2 mm}

$(ii)$  $\mathcal{B}$  is antitall.
\vspace{5 mm}

Proof.} $(i)\Longrightarrow  (ii)$.
Let $Y\subseteq  X$,  $Y \notin \mathcal{B}$.
We take a metric $d$  such that $\mathcal{B}\subseteq  \mathcal{B}_{d}$,  $Y\notin \mathcal{B}_{d}$.
We choose a converging sequence $(z_{n}) _{n\in\omega} $  in $Y$ and put
$Z=\{ z_{n}: n\in \omega\}$.
Then $Z^{\prime}\in \mathcal{B}_{d}$  for each infinite subset $Z^{\prime}$  of $Z$.
Since $\mathcal{B}\subseteq \mathcal{B}_{d}$, we see that $\mathcal{B}$ is antitall.

\vspace{5 mm}

$(ii)\Longrightarrow  (i)$.  Let $Y\subseteq  X$,  $Y\notin \mathcal{B}$.
We choose a countable subset $Z$ of $Y$  such that $Z^{\prime}\notin \mathcal{B}$  for each countable subset $Z^{\prime}$  of $Z$.
Then we endow $X$  with a metric $d$  such that $d(x,y)=1$ for all distinct $x,y\in X\setminus  Z$  and $d$ induces a topology of convergent sequence on $Z$.  Clearly,  $\mathcal{B}\subseteq  \mathcal{B}_{d}$ and $Y\notin \mathcal{B}_{d}$.

\vspace{5 mm}

{\bf Example}.
We take a set $X$  of cardinality $> 2^{\aleph_{0}}$ and denote by $\mathcal{B}$  the  bornology of all finite subsets of $X$.
Evidently,  $\mathcal{B}$  is antitall.
We assume that there is a metric $\rho$  on $X$  such that $\mathcal{B}=\mathcal{B}_{\rho}$.
Since every closed dense subset of $(X,\rho)$ is finite, we conclude that $|X|\leq  2^{\aleph_{0}}$.
Hence, $\mathcal{B}\notin \{ \mathcal{B}_{d}: d$  is a metric on $X\}.$

\section{Tall and antitall  bornologies}

We endow a set $X$  with the discrete topology, identify the Stone-$\check{C}$ech  compactification
 $\beta X$  of $X$ with the set of all ultrafilters on $X$  and denote $X^{\ast}=\beta X \setminus X$,
  so $X^{\ast}$  is the set of
  all free ultrafilters on $X$.
Then the family $\{\bar{A}: A \subseteq X\}$, where $\bar{A}=\{ p\in \beta X: A \in p\} $, forms the base for the topology of $\beta X$.
Given a filter $\varphi$  on  $X$, we denote $\bar{\varphi}= \cap \{ \bar{A} \in \varphi\}$, so
$\varphi$  defines  the closed subset $\bar{\varphi}$ of $\beta X$,  and each non-empty  closed subset $K$  of $\beta X$  can be obtained in this way: $K=\bar{\varphi}$,  $\varphi = \{A \subseteq X:  K\subseteq  \bar{A}\}$.

For an ideal $\mathcal{I}$ in $\mathcal{P}_{X}$, we put
$$\mathcal{I}^{\wedge}=\{p\in\beta G: X\setminus A\in p  \  for  \   each  \   A\in\mathcal{I}\},$$
 and note that $\mathcal{I}$  is a bornology if and only if $\mathcal{I}^{\wedge}\subseteq  X^{\ast}$.

 Using this correspondence between bornologies on $X$  and closed subsets of $X^{\ast}$,  we get
\vspace{5 mm}

{\bf Proposition 1}.  {\it A bornology $\mathcal{B}$ on $X$ is tall if and only if $\mathcal{B}^{\wedge}$ is nowhere dense in $X^{\ast}$. A bornology $\mathcal{B}$ on $X$ is antitall if and only if $\mathcal{B}^{\wedge}$ has a dense subset open in $X^{\ast}$. Every bornology on $X$  is the  intersection of some tall and antitall bornologies}.

\vspace{5 mm}

{\bf Proposition 2}.  {\it For a bornology $\mathcal{B}$ on $X$,  the following statements are equivalent:
\vspace{2 mm}

$(i)$ $\mathcal{B}$ is antitall;
\vspace{2 mm}

$(ii)$ if  $Y\subseteq X$ and $Y\notin \mathcal{B}$ then there exists a function
$f: X\longrightarrow\omega$ such that $f$  is bounded on each member of  $\mathcal{B}$ but $f$ is unbounded on $Y$.}
\vspace{4 mm}

{\bf Remark}. When this note was in the late embryonal state, Taras Banakh noticed that each tall ideal on $X$ does not satisfy  $(ii)$ of Proposition 2. After that, the appearance of antitall ideals was inevitable.
\vspace{4 mm}

Given a bornology  $\mathfrak{I}$ on a set $X$, how to construct the smallest antitall  bornology $\mathcal{B}$ such that $\mathfrak{I}\subseteq\mathcal{B}$?
We make it in topological and functional ways.

Applying Proposition 1, we conclude that  $\mathcal{B} ^{\wedge} = cl \ (int \  \mathfrak{I}^{\wedge})$. So   we take the  filter $\varphi$ such that  $\overline{\varphi}= \mathcal{B} ^{\wedge}$ .
Then  $\mathcal{B}=\{ X\setminus   A:  A  \in \varphi\}$.

By Proposition 2, to get $\mathcal{B}$ it  suffices to join to $\mathfrak{I}$  each subset $B$ of $X$  such that every function  $f: X\longrightarrow\omega $, bounded on every member of $\mathfrak{I}$, is bounded on $B$.

Given an ideal  $\mathfrak{I}$  on $X$,  we can construct some antitall  bornology $\mathcal{B}$  from member of $\mathfrak{I}$ in essentially different way.
 For a family $\mathfrak{F}$ of subsets of  $X$, we denote
$$\mathfrak{F}^{\neg}= \{ A \subseteq X:
\  every \ infinite \ subset  \ of  \  A  \  is  \ not  \ in \  \mathfrak{F} \}. $$
If $\mathfrak{F}$ is inherited by subsets then $\mathfrak{F}^{\neg}$  is a bornology.
 \vspace{5 mm}

{\bf Proposition 3}. {\it  For every ideal $\mathfrak{I}$ on $X$,  $\mathfrak{I}^{\neg}$  is an antitall  bornology.
If $\mathcal{B}$  is an antitall   bornology  then  $(\mathcal{B}^{\neg})^{\neg} = \mathcal{B}$.}
 \vspace{3 mm}

To prove Proposition 3, it suffices to understand a topological sense of the operation $\neg$.
In fact,  $$(\mathfrak{I}^{\neg})^{\wedge} = cl (X^{\ast} \setminus cl \ (int \  \mathfrak{I}^{\wedge})). $$

We conclude with one more topological observation. For a bornology $\mathfrak{I}$  on  $X$,  we endow  $\mathcal{P}_{X}$  with the  topology of uniform
convergence on subsets from $\mathfrak{I}$.
For $Y\in \mathcal{P}_{X}$,  the family
$\{Z\in \mathcal{P}_{X}:  Z\cap A= Y\cap A\}$,
  $A\in \mathfrak{I}$,  is a base of $\tau_{\mathfrak{I}}$ at
  the point $Y$.
We note that $(\mathcal{P}_{X}, \tau_{\mathfrak{I}} )$ is a complete  topological Boolean group
(with the symmetric  difference as  a group operation). We consider $X$ as the subspace $( \{\{ x \}  :  x\in X\}) $   of $(\mathcal{P}_{X}, \tau_{\mathfrak{I}} )$ .
 \vspace{5 mm}

{\bf  Proposition 4.}  {\it For a bornology  $\mathfrak{I}$  on $X$,  the following statements are equivalent:
 \vspace{3 mm}

$(i)$  every infinite subset $Y\subseteq  X$  has an infinite closed discrete subset in $\tau_{\mathfrak{I}}$;
 \vspace{2 mm}

$(ii)$  $\mathfrak{F}_{X}$ is closed (and so  complete) subgroup of $(\mathcal{P}_{X}, \tau_{\mathfrak{I}} ); $
 \vspace{2 mm}

$(iii)$  $\mathfrak{I}$  is tall.}
 \vspace{5 mm}

{\bf  Proposition 5.}  {\it For a bornology  $\mathfrak{I}$  on $X$,  the following statements are equivalent:
 \vspace{3 mm}

$(i)$  each non-closed in $\tau_{\mathfrak{I}}$ subset $Y \subseteq X$  has a convergent (to $\emptyset$ ) sequence;
 \vspace{2 mm}

$(ii)$   $\mathfrak{I}$  is antitall.}

\section{ Classes of bornological spaces}

1.
A set $X$ endowed with a bornology $\mathcal{B}$ is called a {\it bornological space,}   and is denoted by $(X, \mathcal{B})$. Each $A\in \mathcal{B} $  is called   {\it bounded}.

A class of bornological spaces closed under subspaces,
 products and bornologous  images is called a {\it variety.}

A mapping $f: (X,\mathcal{B})\longrightarrow (X^{\prime},\mathcal{B}^{\prime})$
 is {\it bornologous  } if $f(A)\in \mathcal{B}^{\prime}$
  for each  $A\in\mathcal{B}$.

The product of a family of bornological spaces is the Cartesian product of its supports endowed  with the Cartesian product of its bornologies.

We denote by $\mathfrak{M}  _{single}$  the variety of all singletons,  $ \mathfrak{M}  _{bound}$ the variety of all bounded $(\mathcal{B} = \mathcal{P} _{X})$ bornological spaces, and $\mathfrak{M}_{\kappa}$ the variety of all  $\kappa$-bounded bornological spaces. For an infinite cardinal $\kappa$, $(X, \mathcal{B})$  is called
 $\kappa$-{\it bounded}  if    $\mathcal{B}\supseteq [X]^{<\kappa}$.

Applying Theorem 2 from [3] and item 3 below, we
conclude that each variety of bornological spaces  lies in the chain
$$\mathfrak{M}  _{single} \subset \mathfrak{M}  _{bound}\subset \ldots \mathfrak{M}  _{\kappa}\subset\ldots \mathfrak{M}  _{\omega}.$$

We note  that each variety  in this chain,  excepts $ \mathfrak{M}  _{\omega}$, consists of tall bornological spaces.

\vspace{5 mm}

2.   The class of all tall space is closed under subspaces
 and bornologous images   (but not products). The class of all antitall spaces is closed under subspaces,
products ( but not bornologous  images).

For a bornological space $(X, \mathcal{B})$,  we define the {\it hyperbornology} $\exp\mathcal{B}$   on
$\mathcal{B}$: the family $\{Y\in\mathcal{B}: Y\subseteq A\}$, $A\in \mathcal{B}$
 is the base of $\exp\mathcal{B}$.

We show that if $\mathcal{B}$  is antitall then  $\exp\mathcal{B}$  is antitall.
 Let $\mathcal{A}$  be a family of bounded subsets of  $(X, \mathcal{B})$
 such that $\mathcal{A}\notin\exp\mathcal{B}$.
 We   denote $Y=\cup\mathcal{A}$.
  Since  $Y\notin\mathcal{B}$
  and $\mathcal{B}$ is antitall, there is a countable $Z \subseteq Y$  such that each
  countable  subset of $Z$  is unbounded in
  $(X,\mathcal{B})$.
For each $z\in  Z$, we pick  $A_{z}\in  \mathcal{A}$  such that  $z\in  A _{z}$.
Then $\{A_{z}: z\in Z\}$
witness that $\exp\mathcal{B}$  is antitall.

On the tall hand, let $\mathcal{B}$ be a bornology on a countable set $X$  such that $\mathcal{B}\neq\mathcal{P}_{X}$.
We write $X$ as the union of increasing chain $\{X_{n}: n\in\omega\}$ of finite subsets. Clearly, only finite subsets of
$\{X_{n}: n\in\omega\}$ are in $\exp\mathcal{B}$.
Hence, $\exp\mathcal{B}$ is not tall.

\vspace{5 mm}

3. Following [5], we say that a family $\mathcal{E}$  of subsets of $X \times X$
is a coarse structure on a set $X$  if

\begin{itemize}
\item{}
each  $\varepsilon\in\mathcal{E}$ contains the  diagonal $\triangle_{X}$, $\triangle_{X}=\{(x,x): x\in X\}$;

\item{}  if $\varepsilon,\delta\in\mathcal{E}$ then $\varepsilon\circ\delta\in\mathcal{E}$ and $\varepsilon^{-1} \in\mathcal{E}$ where $\varepsilon\circ\delta= \{(x,y): \exists z ((x,z)\in\varepsilon, (z,y)\in\delta )\}$,
    $\varepsilon^{-1}= \{(y,x): (x,y)\in\varepsilon\}$;

\item{}  if $\varepsilon\in\mathcal{E}$ and $\triangle_{X}\subseteq \varepsilon^{\prime}\subseteq \mathcal{E}$  then
$\varepsilon^{\prime}\in\mathcal{E}$;

\item{} for any $x,y\in X$, there is $\varepsilon\in\mathcal{E}$ such that $(x,y)\in\varepsilon$.

\end{itemize}

The pair  $(X, \mathcal{E})$ is called a {\it coarse space}. For $x\in X$ and $\varepsilon\in\mathcal{E}$,
we denote $B(x,\varepsilon)=\{y\in X: (x, y)\in\varepsilon\} $
 and say that $B(x,\varepsilon)$
 is a  {\it ball of  radius $\varepsilon$  around $x$}.
We note that a coarse space can be considered as an asymptotic counterpart of a uniform
spaces and could be defined in terms of balls, see [4]. In this case a coarse space is called a {\it ballean}.

Let $(X, \mathcal{E})$ be a coarse space.
 A subset $Y$ of  $X$  is called {\it bounded} if there exist  $x\in  X$  and $\varepsilon\in\mathcal{E}$
  such that  $Y\subseteq  \mathcal{B}(x, \mathcal{E})$.
The family of all bounded  subsets of  $(X, \mathcal{E})$  is a bornology.
On the other side, for every bornology $\mathcal{B}$
 on   $X$,  there is the smallest (by inclusion) coarse structure $\mathcal{E}_{\mathcal{B}}$
  such that $\mathcal{B}$ is a bornology  of all bounded subsets of   $(X, \mathcal{E}_{\mathcal{B}})$.
A coarse space  $(X,\mathcal{E})$  is of the form $(X, \mathcal{E}_{\mathcal{B}})$
 if and only if $(X,\mathcal{E})$
   is thin:  for every $\varepsilon\in\mathcal{E}$,
    there is a bounded subset  $A$  of $(X,\mathcal{E})$  such that $B(x,\varepsilon)=\{x\}$
      for all $x\in X \backslash  A$.

\bibliography{mybibfile}

\vspace{3 mm}






\begin{thebibliography}{2}


\bibitem{b1}{ R.Engelking}, {\it General Topology, } 2nd edition, PWN, Warszawa, 1985.


\bibitem{b2}{ H. Hogbe-Nlend,} {\it Les racines historiques de la bornologie moderne}, Seminare Choquet,
{\bf 10} (1970-71), no1,1-7.

\bibitem{b3}{ I. Protasov,} {\em Varieties of coarse spaces},  Axioms, 2018,  7, 32.

\bibitem{b4}{ I. Protasov, M. Zarichnyi,} {\em General Asymptology},   Math. Stud. Monogr. Ser., Vol. 12, VNTL, Lviv, 2007,  pp. 219.

\bibitem{b5}{ J. Roe, }{\it Lectures on Coarse Geometry}, AMS University Lecture Ser {\bf 31}, Providence, R.I., 2003, pp. 176.



\end{thebibliography}

\vspace{6 mm}

CONTACT INFORMATION

I.~Protasov: \\
Faculty of Computer Science and Cybernetics  \\
        Kyiv University  \\
         Academic Glushkov pr. 4d  \\
         03680 Kyiv, Ukraine \\ i.v.protasov@gmail.com

\end{document}